\documentclass[12pt]{amsart}
\usepackage[utf8]{inputenc}
\usepackage{amsmath,amssymb,amsthm}
\usepackage[all]{xy}
\usepackage{xcolor}
\usepackage{tikz,tikz-cd}
\usepackage{lipsum} 
\usepackage{mathtools}

\DeclareMathOperator{\Ex}{\mathrm{Ex}}
\def\M{\mathbf M}
\def\X{\Tilde X}
\def\B{\Tilde B}
\def\A{\Tilde A}
\def\D{\Tilde D}
 
\newtheorem{prop}[equation]{Proposition} 
\newtheorem{theorem}[equation]{Theorem}

\newtheorem{lemma}[equation]{Lemma}

\newtheorem{corollary}[equation]{Corollary}
\theoremstyle{definition}
\newtheorem{definition}{Definition}[section]
\theoremstyle{remark}

\title{Flops and minimal models of generalized pairs}
\author{Priyankur Chaudhuri}
\address{School of Mathematics, Tata Institute of Fundamental Research, Homi Bhabha Road, Colaba,
Mumbai 400005}
\email{pkurisibang@gmail.com}

\begin{document}
\begin{abstract}
    We show that given any two minimal models of a generalized lc pair, there exist small birational models which are connected by a sequence of symmetric flops. We also present some applications of this result. 
\end{abstract}

\maketitle

\section{Introduction}

In \cite{Ka}, Kawamata showed that any two $\mathbb{Q}$-factorial terminal pairs $(X,B)$ and $(X^{'},B^{'})$ with $K_X+B$ and $K_{X^{'}}+B^{'}$ nef which are birational to each other can be connected by a sequence of small birational maps known as flops. Hashizume \cite{Has} proved a version of this result for log canonical pairs that are not necessarily $\mathbb{Q}$-factorial. In view of recent developments in the minimal model program of generalized pairs, it is natural to ask if similar results hold for minimal models of generalized pairs. In this article, we closely follow the ideas of Hashizume \cite{Has} and use some recent results on generalized pairs due to Hacon, Liu and Xie (\cite{HL}, \cite{LX}) to extend this result to generalized log canonical pairs.

\begin{theorem}
    
    Suppose $(X,B+\mathbf{M})/S$ and $(X^{'}, B^{'}+\mathbf{M})/S$ are two generalized log canonical pairs such that $K_X+B+\mathbf{M}_X$ and $K_{X^{'}}+B^{'}+\mathbf{M}_{X^{'}}$ are nef over $S$, $\mathbf{M}_X$ and $\mathbf{M}_{X^{'}}$ are $\mathbb{R}$-Cartier and there exists a small birational map $\phi:X \dashrightarrow X^{'}$ over $S$ such that
    
    \begin{itemize}
        \item $B^{'}=\phi_*B$ and $\mathbf{M}_{X^{'}}= \phi_*\mathbf{M}_X$,
        \item there exists $ U\subset X$ open such that $\phi|_U$ is an isomorphism and all glc centers of $(X,B+\mathbf{M})$ intersect $U$
    \end{itemize}
    then (possibly after exchanging $X$ and $X^{'}$), there exist small birational morphisms from normal quasi-projective varieties $(\Tilde{X}, \Tilde{B}+\mathbf{M}) \xrightarrow{\Tilde{\psi}}(X,B+\mathbf{M})$ and $(\Tilde{X^{'}}, \Tilde{B^{'}}+\mathbf{M}) \xrightarrow{\Tilde{\psi^{'}}}(X^{'}, B^{'}+\mathbf{M})$ such that the induced birational map $(\Tilde{X}, \Tilde{B}+\mathbf{M}) \dashrightarrow  (\Tilde{X^{'}}, \Tilde{B^{'}}+\mathbf{M}) $ can be written as a composition of a finite sequence of symmetric flops (see definition \ref{deflop}) over $S$ with respect to $K_{\Tilde{X}}+\Tilde{B}+\mathbf{M}$. 

\end{theorem}

Note that the open subset $U \subset X$ in the theorem exists if $(X,B+\M)$ and $(X^{'}, B^{'}+\M)$ arise as outputs of two MMP's on a generalized log canonical pair (see Lemma \ref{mm}). \\

The outline of the proof is as follows: let $A^{'}$ be a general ample divisor on $X^{'}$ such that $(X^{'},B^{'}+A^{'}+ \mathbf{M})$ is generalized lc. Let $A$ denote the strict transform of $A^{'}$ on $X$. In the $\mathbb{Q}$-factorial generalized klt case, we can clearly choose $A^{'}$ such that $(X^{'}, B^{'}+A^{'}+\mathbf{M})$ and $(X, B+A+\mathbf{M})$ are both generalized klt. $(X^{'},B^{'}+A^{'}+\mathbf{M})$ is then the log canonical model of $(X,B+A+\mathbf{M})$ and one can show by the arguments of Kawamata \cite{Ka} that there exists a sequence of symmetric flops connecting them. If we drop $\mathbb{Q}$-factoriality, then the flops only take us to a good minimal model of $(X,B+A+\mathbf{M})$ which is then a small birational model of the log canonical model $(X^{'}, B^{'}+A^{'}+\mathbf{M})$.\\

In the generalized log canonical case, it is not clear if we can choose an ample $A^{'}$ such that $(X^{'},B^{'}+A^{'}+\mathbf{M})$ and $(X,B+A+\mathbf{M})$ are both generalized lc. However if we take a generalized dlt modification $(\hat{X}, \hat{B}+ \mathbf{M}) \xrightarrow{\mu}(X,B+\mathbf{M})$, if $\hat{A}$ denotes the proper transform of $A^{'}$ on $\hat{X}$ then $(\hat{X}, \hat{B}+ \hat{A}+ \mathbf{M})$ be made generalized log canonical. Using a recent result of Liu and Xie \cite{LX}, we next show that $(\hat{X}, \hat{B}+ \hat{A}+ \mathbf{M})$ has a generalized log canonical model $(\Tilde{X}, \Tilde{B}+\Tilde{A}+\mathbf{M})$ over $X$. Then the induced birational morphism $\Tilde{\mu}: \Tilde{X} \rightarrow X$ is small and we can run a sequence of flops on $\Tilde{X}$ as above. Note that these flops don't necessarily preserve Picard rank. \\

We now discuss some applications of this result. As noted by Hashizume \cite[Remark 4.7]{Has}, flips for log canonical pairs don't necessarily preserve the cohomology groups of the structure sheaf. However, these cohomology groups do agree for all minimal models obtained from a given log canonical pair by running various MMP's \cite[Theorem 1.2]{Has}. We show that this continues to hold in the setting of generalized pairs. Our other application is concerned with the invariance of Cartier index. In general, the log canonical divisors of two minimal models of a given lc pair need not have the same Cartier index \cite[Example 4.8]{Has}. However, if two minimal models arise by running two MMP's on a given lc pair, then their log canonical divisors have the same Cartier index \cite[Theorem 1.2]{Has}. This holds for generalized pairs as well. More generally, we have:

\begin{theorem}
    Suppose $(X,B+\mathbf{M})/S$ and $(X^{'}, B^{'}+\mathbf{M})/S$ are two generalized log canonical pairs with structure morphisms $\pi: X \to S$ and $\pi^{'}: X^{'}\to S$ and such that $K_X+B+\mathbf{M}_X$ and $K_{X^{'}}+B^{'}+\mathbf{M}_{X^{'}}$ are nef over $S$, $\mathbf{M}_X$ and $\mathbf{M}_{X^{'}}$ are $\mathbb{R}$-Cartier and there exists a small birational map $\phi:X \dashrightarrow X^{'}$ over $S$ such that
    
    \begin{itemize}
        \item $B^{'}=\phi_*B$ and $\mathbf{M}_{X^{'}}= \phi_*\mathbf{M}_X$,
        \item there exists $ U\subset X$ open such that $\phi|_U$ is an isomorphism and all glc centers of $(X,B+\mathbf{M})$ intersect $U$.
    \end{itemize}
    Then we have the following:

    \begin{enumerate}
        \item $R^p \pi_* \mathcal{O}_X \cong R^p \pi^{'}_*\mathcal{O}_{X^{'}}$ for all $p >0$. In particular, if $S$ is a point, then $H^i(X, \mathcal{O}_X) \cong H^i(X^{'}, \mathcal{O}_{X^{'}})$ for all $i>0$,
        \item $K_X+B+\M_X$ and $K_{X^{'}}+B^{'}+\M_{X^{'}}$ have the same Cartier index.
        \end{enumerate}

\end{theorem}
\section{Preliminaries}

\begin{definition}
(Generalized pairs and their singularities \cite[Definition 1.4, 4.1]{BZ}) A \emph{generalized sub-pair} $(X,B+\mathbf{M})/S$ consists of a normal quasi-projective variety $X$ equipped with a projective morphism to a variety $S$, an $\mathbb{R}$-divisor $B$ and an $\mathbb{R}$-b-divisor $\mathbf{M}$ on $X$ such that:
\begin{itemize}
\item $K_X+B+\mathbf{M}_X$ is $\mathbb{R}$-Cartier.
\item $\mathbf{M}$ is b-nef NQC i.e. it descends to a nef $\mathbb{R}$-divisor on a birational model of $X$ where it can be written as a real linear combination of nef $\mathbb{Q}$-Cartier divisors.
\end{itemize}
When $B \geq 0$, we drop the prefix sub.\\

For any prime divisor $E$ and an $\mathbb{R}$-divisor $D$ on $X$, let mult$_E(D)$ denote the multiplicity of $E$ along $D$. \\

Let $(X,B+\mathbf{M})/S$ be a generalized sub-pair and $Y \xrightarrow{\mu} X$ a log resolution of $(X,B)$ such that $\mathbf{M}$ descends to $Y$. Let $B_Y$ be defined by $K_Y+B_Y+\mathbf{M}_Y = \mu^*(K_X+B+\mathbf{M}_X)$. We say that $(X,B+\mathbf{M})$ is \emph{generalized sub-klt} (resp. \emph{generalized sub-lc}) if every coefficient of $B_Y$ is less than $1$ (resp. $\leq 1$). \newline 
The \emph{discrepancy} of a prime divisor $D$ on $Y$ with respect to $(X,B+\mathbf{M})$ is defined and denoted by $a(D,X,B+\mathbf{M}):= -$mult$_D(B_Y)$. Thus $(X, B+\mathbf{M})$ is generalized sub-klt (resp. sub-lc) if $a(D,X,B+\mathbf{M}) > -1$ (resp. $\geq -1$) for every log resolution $Y$ as above and every prime divisor $D$ on $Y$. The discrepancy divisor is defined by $A_Y(X,B+\mathbf{M})= \Sigma_D a(D,X,B+\mathbf{M})D$.\\

We say that a generalized sub-lc pair $(X, B+\mathbf{M})$ is \emph{generalized sub-dlt} if there exists a closed subset $V \subset X$ such that 
\begin{itemize}
    \item $X\setminus V$ is smooth and $B|_{X\setminus V}$ is simple normal crossing,
    \item for any prime divisor $E$ over $X$ such that $a(E, X,B+\mathbf{M})=-1$, we have center$_XE \not \subset V$ and $\mathbf{M}$ descends to a nef divisor on center$_X(E\setminus V)$.
\end{itemize}
In case $\mathbf{M}=0$, and $(X,B+\mathbf{M})$ is generalized sub-lc (resp. generalized sub-klt etc), then we say $(X,B)$ is sub-lc (resp. sub-klt etc).\\

From now on, we use \emph{glc, gklt} and \emph{gdlt} to denote generalized lc, generalized klt and generalized dlt respectively.\\

Suppose $(X, B+\mathbf{M})$ is glc. A \emph{glc place} is a prime divisor $E$ over $X$ such that $a(E, X, B+ \mathbf{M}) =-1$. A \emph{glc center} of $(X,B+\mathbf{M})$ is the center on $X$ of a glc place of $(X,B+\mathbf{M})$.
\end{definition}

We will need the following simple consequence of the negativity lemma:

\begin{lemma}\label{M}
    Let $(X, B+\M)$ be a generalized lc pair with $\M_X$ $\mathbb{R}$-Cartier. Then $(X,B)$ is lc.

    \begin{proof}
        Let $\pi:Y \to X$ be a log resolution of $(X,B)$ such that $\M_Y$ descends to a nef divisor on $Y$. Let $B_Y$ be the divisor on $Y$ defined by $K_Y+B_Y= \pi^*(K_X+B)$. By negativity lemma, we have $\M_Y= \pi^*\M_X-E$ for some effective $\pi$-exceptional divisor $E$. Then  we have $K_Y+B_Y+E+\M_Y = \pi^*(K_X+B+\M_X)$. Since $(X, B+\M)$ is glc, the coefficients of $B_Y+E$ are atmost $1$ and hence all coefficients of $B_Y$ are also atmost $1$. Thus $(X,B)$ is lc. 
    \end{proof}
\end{lemma}

\begin{definition}
(Generalized models \cite[Definition 2.1]{Bi}, \cite[Definition 2.21]{HL}) Let $(X, B+\mathbf{M})/S$ be a generalized log canonical pair. A generalized pair $(X^{'}, B^{'}+\mathbf{M})/S$ equipped with a birational map $\phi: X \dashrightarrow X^{'}$ over $S$ is called a \emph{generalized weak log canonical model} of $(X,B+\mathbf{M})/S$ if 
\begin{itemize}
    \item $B^{'}=\phi_*(B)+E$, where $E$ is the reduced $\phi^{-1}$-exceptional divisor,
    \item $K_{X^{'}}+B^{'}+\mathbf{M}_{X^{'}}$ is nef over $S$,
    \item $(X^{'},B^{'}+\mathbf{M})$ is generalized lc,
    \item For any prime divisor $D$ on $X$ which is exceptional over $X^{'}$, we have $a(D,X, B+\mathbf{M}) \leq a(D,X^{'}, B^{'}+ \mathbf{M})$.
\end{itemize}

A generalized weak log canonical model $(X^{'},B^{'}+\mathbf{M})$ of $(X,B+\mathbf{M})/S$ is called a \emph{generalized log canonical model} (glc model in short) of $(X,B+\mathbf{M})$ if $K_{X^{'}}+B^{'}+\mathbf{M}_{X^{'}}$ is ample over $S$.\\
 A generalized weak log canonical model $(X^{'},B^{'}+\mathbf{M})/S$ as above is called a \emph{generalized minimal model} of $(X,B+\mathbf{M})/S$ if for any prime divisor $D$ on $X$ which is exceptional over $X^{'}$, we have $a(D,X, B+\mathbf{M}) < a(D,X^{'}, B^{'}+ \mathbf{M})$.\\
 If $(X^{'},B^{'}+\mathbf{M})/S$ is a generalized minimal model of $(X,B+\mathbf{M})/S$ such that $K_{X^{'}}+B^{'}+\mathbf{M}_{X^{'}}$ is semiample over $S$, then we say that $(X^{'}, B^{'}+ \mathbf{M})/S$ is a \emph{generalized good minimal model} for $(X,B+\mathbf{M})/S$. 
 \newline Since we will mainly deal with generalized pairs throughout this article, we will drop the term 'generalized' and just use phrases like minimal model etc.

\end{definition}
Now we recall the definitions of flips and flops for generalized pairs. 

\begin{definition}
($D$-flips \cite[Definition 8.8]{Bi2}) Let $X$  be a normal  variety equipped with a projective morphism $X \rightarrow S $. Let $D$ be an $\mathbb{R}$-Cartier divisor on $X$. A birational morphism $f:X \rightarrow V$ over $S$ where $V$ is projective over $S$ is called a \emph{D-flipping contraction} over $S$ if $\rho(X/V)=1$, $f$ is small (i.e. has exceptional locus of codim at least $2$) and $-D$ is ample over $V$. Let $f^{'}: X^{'} \rightarrow V$ be a projective birational morphism over $S$ from a normal variety $X^{'}$ and $\phi: X \dashrightarrow X^{'}$ the induced birational map. Then $f^{'}$ is a \emph{flip} of $f$ if $f^{'}$ is small and $\phi_*D$ is $\mathbb{R}$-Cartier  and ample over $V$.

 \end{definition}

 \begin{definition}(Extremal contractions \cite[Definition 3.34]{KM})
     A contraction morphism $X \rightarrow Y$ is called \emph{extremal} if for any two Cartier divisors $D_1$ and $D_2$ on $X$, there exist integers $a,b$ not both zero such that $aD_1-bD_2$ is linearly equivalent to the pullback of some Cartier divisor on $Y$.
 \end{definition}
\begin{definition}
    
(flops for generalized pairs)\label{deflop} Let $(X, B+\mathbf{M})/S$ be a generalized lc pair where $X$ is projective over $S$. A \emph{flop} for $K_X+B+\mathbf{M}$ over $S$ consists of the diagram 

\begin{equation*}
\xymatrix@R=16pt
{
X\ar@{-->}^{\varphi}[rr]\ar[dr]^{f}\ar[ddr]&&X' \ar[dl]_{f'}\ar[ddl]\\
&V\ar[d]\\
&S
}
\end{equation*}
where $f$ is a $D$-flipping contraction over $S$, $f^{'}$ is a flip of $f$ and $K_X+B+\mathbf{M} \equiv _{V}0$. The flop is called a \emph{symmetric flop} \cite[Definition 2.4]{Has} if both $f$ and $f^{'}$ are extremal contractions.
\end{definition}

We have the following criterion for a $K_X+B+\mathbf{M}$-flop to be symmetric.

\begin{lemma}\cite[Lemma 2.6]{Has}\label{sym} Let 

\begin{equation*}
\xymatrix@R=16pt
{
X\ar@{-->}^{\varphi}[rr]\ar[dr]^{f}\ar[ddr]&&X' \ar[dl]_{f'}\ar[ddl]\\
&V\ar[d]\\
&S
}
\end{equation*}
be a $K_X+B+\mathbf{M}_X$-flop over $S$ (where $f$ is a $D$-flipping contraction and $f^{'}$ its flip) such that $\mathbf{M}_X$ is $\mathbb{R}$-Cartier and 

\begin{enumerate}
    \item $\rho (X/S)= \rho(X^{'}/S)$,
    \item There exists an effective $\mathbb{R}$-Cartier divisor $E$ on $X$ such that $(X,B+E+\mathbf{M})$ is glc and $-(K_X+B+E+\mathbf{M}_X)$ is ample over $V$.
\end{enumerate}
Then the flop is symmetric (see Definition \ref{sym}). Furthermore, there exists an effective $\mathbb{R}$-Cartier divisor $F^{'}$ on $X^{'}$ such that $(X^{'},\phi_*B+F^{'}+\mathbf{M})$ is glc and $-(K_X+\phi_*B+F^{'}+\mathbf{M}_{X^{'}})$ is ample over $V$.

  \begin{proof}
      Since $(X,B+E+\mathbf{M})$ is glc and $-(K_X+B+E+\mathbf{M}_X)$ is ample over $V$, by taking $H$ as a general member of the $\mathbb{R}$-linear system of $-(K_X+B+E+\mathbf{M}_X)$ over $V$, we can make $(X,B+E+H+\mathbf{M})$ glc and $K_X+B+E+H+\mathbf{M}_X \equiv _{V} 0$. Since $\mathbf{M}_X$ is $\mathbb{R}$-Cartier, by contraction theorem for glc pairs (\cite[Theorem 1.3 (4)]{HL}), we conclude that $K_X+B+E+H+\mathbf{M}_X \sim _{\mathbb{R},V} 0$. Let $G$ be any $\mathbb{R}$-Cartier divisor on $X$. Since $\rho (X/V)=1$, there exists $r \in \mathbb{R}$ such that $G-rD \equiv _{V} 0$ and then $G-rD \sim_{\mathbb{R}, V}0$ as above. From this, it easily follows that $f$ is an extremal contraction. To show the same for $f^{'}$, we will need to produce a glc structure $(X^{'},\Delta^{'}+\mathbf{M})$ on $X^{'}$ such that $-(K_{X^{'}}+\Delta^{'}+\mathbf{M}_{X^{'}})$ is ample over $V$. By taking $\phi_*$, we get $\phi_*G-r \phi_*D \sim _{\mathbb{R},V} 0$. Since $\phi_*D$ is $\mathbb{R}$-Cartier, so is $\phi_*G$. This gives $\phi_* :N^1(X/S) \rightarrow N^1(X^{'}/S)$. Since $\phi$ is small, $\phi_*$ is injective. (Indeed, let $p:\Tilde{X} \rightarrow X$ and $q: \Tilde{X} \rightarrow X^{'}$ resolve $\phi$. Let $G$ be an $\mathbb{R}$-Cartier divisor on $X$ such that $\phi_*G \equiv_S 0$. Write $p^*G=q^*\phi_*G+E$ where $E$ is exceptional over both $X$ and $X^{'}$. Then $q^*\phi_*G \equiv _{S} 0$ and thus $G=p_*q^*\phi_*G \equiv _S 0$). Since $\rho(X/S)=\rho(X^{'}/S)$, it follows that $\phi_*$ is actually an isomorphism. \\

      By the definition of $H$, it follows that $K_{X^{'}}+\phi_*(B+E+H)+\mathbf{M}_{X^{'}}$ is $\mathbb{R}$-Cartier. Also $K_{X^{'}}+\phi_*(B+E+H)+\mathbf{M}_{X^{'}} \sim_{\mathbb{R},V} 0$ and $(X^{'}, \phi_*(B+E+H)+\mathbf{M})$ is glc. Since $-E$ is ample over $V$ (follows from the fact that $K_X+B+\M_X$ is trivial over $V$ and $-(K_X+B+E+\M_X)$ is ample over $V$), it follows that $E \sim_{\mathbb{R},V} \alpha D$ for some $ \alpha >0$. Then $\phi_*E \sim _{\mathbb{R},V}\alpha \phi_*D$. Since $\phi_*D$ is ample over $V$ ($f^{'}$ is a $D$-flip), so is $\phi_*E$. Thus 
      \begin{center}
          $-(K_{X^{'}}+ \phi_*B+\phi_*H+\mathbf{M}_{X^{'}})= -(K_{X^{'}}+\phi_*(B+E+H)+\mathbf{M}_{X^{'}})+\phi_*E$
      \end{center}
      is ample over $V$. Clearly $(X^{'}, \phi_*(B+H)+\mathbf{M})$ is glc since $(X^{'},\phi_*(B+E+H)+\mathbf{M})$ is glc and $\phi_*E \geq 0$. Then we can take $F^{'}:= \phi_*H$.\\

  \end{proof}
\end{lemma}

\section{Relation between minimal models}
\begin{lemma} \label{mm}
    Let $(X,B+\mathbf{M})/S$ be a glc pair. Let $\phi_1:(X,B+\mathbf{M})\dashrightarrow (Y,B_Y+\mathbf{M})$ and $\phi_2:(X,B+\mathbf{M}) \dashrightarrow (Y^{'},B_{Y^{'}}+\mathbf{M})$ be two minimal models of $(X,B+\mathbf{M})/S$. Let $\phi: Y \dashrightarrow Y^{'}$ be the induced birational map. Then we have the following:

    \begin{itemize}
        \item If $(Y,B_Y+\mathbf{M})$ is a good minimal model of $(X,B+\mathbf{M})/S$ then so is $(Y^{'},B_{Y^{'}}+\mathbf{M})$
        
        \item In case $(Y,B_Y+\mathbf{M})$ and $(Y^{'},B_{Y^{'}}+\mathbf{M}_{Y^{'}})$ are obtained by a sequence of steps of a $(K_X+B+\mathbf{M}_X)$-MMP over $S$, then $\phi_*B_Y=B_{Y^{'}}$, $\phi_*\mathbf{M}_Y= \mathbf{M}_{Y^{'}}$. Moreover, there exists $U \subset Y$ open such that $\phi|_{U}$ is an isomorphism and all glc centers of $(Y,B_Y+\mathbf{M}_Y)$ intersect $U$.
    \end{itemize}
   
    \begin{proof}
        Let $W$ be a smooth resolution of indeterminacy of $\phi_1$ and $\phi_2$ with induced morphisms $p:W \rightarrow Y$, $q:W\rightarrow X$ and $r:W \rightarrow Y^{'}$ such that $\mathbf{M}$ descends to $W$. Let 
        \begin{center}
            $q^*(K_X+B+\mathbf{M}_X)=r^*(K_{Y^{'}}+B_{Y^{'}}+\mathbf{M}_{Y^{'}})+E^{'} \vspace{2 mm}
            \newline =p^*(K_Y+B_Y+\mathbf{M}_Y)+E$.
        \end{center}
        Then \begin{center}
        $E=A_W(Y,B_Y+\mathbf{M})-A_W(X,B+\mathbf{M})$ and
       \vspace{3mm}
         $E^{'}=A_W(Y^{'},B_{Y^{'}}+\mathbf{M})-A_W(X,B+\mathbf{M})$ ($*$).
        \end{center}
        We claim that $E \geq 0$ and is exceptional over $Y$. Indeed, let $D$ be a component of $E$. Then by ($*$), mult$_DE = a(D,Y,B_Y+\mathbf{M}_Y)-a(D,X,B+\mathbf{M})$.\\
        
        Suppose $D$ is not exceptional over $X$. If it is not exceptional over $Y$ either, then $a(D,Y,B_Y+\mathbf{M})=a(D,X,B+\mathbf{M})$, thus forcing mult$_DE=0$ which is a contradiction. Thus we may assume $D$ is exceptional over $Y$. Then by the definition of minimal models, $a(D,Y,B_Y+\mathbf{M})>a(D,X,B+\mathbf{M})$, so mult$_DE>0$. We conclude that $q_*E \geq 0$ (because components $D$ which are exceptional over $X$ map to $0$ and those that are not have positive coefficient in $E$). Now 
        \begin{center}
            $-E=p^*(K_Y+B_Y+\mathbf{M}_Y)-q^*(K_X+B+\mathbf{M}_X)$ \end{center}
            is nef over $X$. Thus $E \geq 0$ by negativity lemma. Now we show that $E$ is exceptional over $Y$. Suppose there exists a component $D$ of $E$ not exceptional over $Y$. If $D$ is not exceptional over $X$, then mult$_DE=0$ as above which is impossible. Thus $D$ is exceptional over $X$ which means $p(D)$ is $\phi_1^{-1}$-exceptional. Since $B_Y=\phi_{1*}B+\Ex(\phi_1^{-1})_{red}$, it follows that $a(D,Y,B_Y+\mathbf{M}_Y)=-1$. Since $E\geq 0$, by ($*$) above and the fact that $(X,B+\mathbf{M})$ is glc, it follows that $a(D,X,B+\mathbf{M})=-1$. This again forces mult$_DE=0$ which is impossible. Thus $E$ is exceptional over $Y$. This proves our claim. \\

            Similar arguments show that $E^{'} \geq 0$ and is exceptional over $Y^{'}$. Thus we have $r_*(E-E^{'}) \geq 0$. Since $E^{'}-E =p^*(K_Y+B_Y+\mathbf{M}_Y)-r^*(K_{Y^{'}}+B_{Y^{'}}+\mathbf{M}_{Y^{'}})$, $E^{'}-E$ is nef over $Y^{'}$. Then $E \geq E^{'} $ by negativity lemma. Similarly, we can show that $E^{'} \geq E$. Thus, we conclude that $E^{'}=E$. This shows that $p^*(K_Y+B_Y+\mathbf{M}_Y)=r^*(K_{Y^{'}}+B_{Y^{'}}+\mathbf{M}_{Y^{'}})$. Thus one of the minimal models is good iff the other is good.\\

            Now we prove the second assertion of the lemma. Let $c_Y(E)$ be a glc center of $(Y,B_Y+\mathbf{M})$ and let $W$ be a common birational model of $X$, $Y$ and $Y^{'}$ as above such that $E$ is a prime divisor on $W$. Then $a(E,Y,B_Y+\mathbf{M}) =-1$. Since $a(E,X,B+\mathbf{M}) \leq a(E,Y,B_Y+\mathbf{M})$ and $a(Y,B_Y+\mathbf{M}) =a(Y^{'},B_{Y^{'}}+\mathbf{M})$ as we had observed above, it follows that $-1=a(E,X,B+\mathbf{M})=a(E,Y^{'},B_{Y^{'}}+\mathbf{M})$ as well. Now by construction of the $(K_X+B+\mathbf{M}_X)$-MMP (see proof of \cite[Lemma 3.38]{KM} for details of the arguments), the discrepancy of a prime divisor $E$ over $X$ strictly increases iff $c_X(E)$ is contained in the non-isomorphic locus of the MMP.\\

            Thus if $V \subset X$ (resp. $V^{'} \subset X$) is the largest open subset on which $\phi_1$ (resp. $\phi_2$) is an isomorphism, it follows that $c_X(E) \cap V \neq \emptyset $ and $c_X(E) \cap V^{'} \neq \emptyset$ (for all glc centers $c_X(E)$). Since $c_X(E)$ is connected, $c_X(E) \cap V \cap V^{'} \neq \emptyset $. So we can take $U = \phi_1(V \cap V^{'})$.\\

            Since both $\phi_1$ and $\phi_2$ are birational contractions in this case, it follows that $\phi_*B_Y=B_{Y^{'}}$ and $\phi_*\mathbf{M}_Y=\mathbf{M}_{Y^{'}}$.

    \end{proof}
\end{lemma}

\begin{prop} \cite[Proposition 3.1]{Has} \label{tech}
    Let $(X,B+\mathbf{M})/S$ and $(X^{'},B^{'}+\mathbf{M})/S$ be two glc pairs such that $K_X+B+\mathbf{M}_X$ and $K_{X^{'}}+B^{'}+\mathbf{M}_{X^{'}}$ are nef over $S$. Suppose $\phi:X \dashrightarrow X^{'}$ is a small birational map over $S$ such that $B^{'}= \phi_*B$, $\mathbf{M}_{X^{'}}= \phi_*\mathbf{M}_X$ and there is an open subset $U \subset X$ such that $\phi$ is an isomorphism on $U$ and all glc centers of $(X,B+\mathbf{M})$ intersect $U$. \newline
    Then there exists a small projective morphism $\psi:\Tilde{X} \rightarrow X$ from a normal quasi-projective variety such that \begin{itemize}
        \item $\psi$ is an isomorphism over $U$,
        \item there is an ample $\mathbb{R}$-divisor $A^{'} \geq 0$ on $X^{'}$ such that $(X^{'},B^{'}+A^{'}+\mathbf{M})$ is glc and if $\Tilde{A}$ is the birational transform of $A^{'}$ on $\Tilde{X}$, then $K_{\Tilde{X}}+\psi_*^{-1}B+\Tilde{A}+\mathbf{M}$ is $\mathbb{R}$-Cartier and $(\Tilde{X}, \psi_*^{-1}B+\Tilde{A}+\mathbf{M})$ is glc.
    \end{itemize}

    \begin{proof}
        Let $Y \xrightarrow{f} X$, $Y \xrightarrow{g}X^{'}$ denote a log resolution of $(X,B)$ that resolves $\phi$. Let $ \Gamma := f_*^{-1}B+\Ex(f)_{red}$. Let $A^{'}$ be an ample $\mathbb{R}$-divisor on $X^{'}$ such that $(X^{'}, B^{'}+A^{'}+\mathbf{M})$ and $(Y, \Gamma+g^*A^{'}+\mathbf{M})$ are both glc. By running a $(K_Y+\Gamma+\mathbf{M})$-MMP over $X$ with scaling of an ample divisor, we construct a $\mathbb{Q}$-factorial g-dlt modification $(Y^{'},\Gamma^{'}+\mathbf{M})$ of $(X,B+\mathbf{M})$ (see the proof of \cite[Lemma 4.5]{BZ} for details). Let $A_{Y^{'}}$ be the birational transform of $g^*A^{'}$ on $Y^{'}$. The map $Y \dashrightarrow Y^{'}$ is also a sequence of steps of a $(K_Y+\Gamma+t_0g^*A^{'}+\mathbf{M})$-MMP over $X$ for $0<t_0 \ll 1$. In particular, this implies that there exists $t \in (0,t_0)$ such that $(Y^{'}, \Gamma^{'}+tA_{Y^{'}}+\mathbf{M})$ is glc and all glc centers of $(Y^{'},\Gamma^{'}+tA_{Y^{'}}+\mathbf{M})$ are also glc centers of $(Y^{'},\Gamma^{'}+\mathbf{M})$ (since running an MMP does not create new glc centers). In particular, all glc centers of $(Y^{'},\Gamma^{'}+tA_{Y^{'}}+\mathbf{M})$ intersect $f^{'-1}(U)$, where $f^{'}:Y^{'} \rightarrow X$ is the induced morphism.\\

        We claim that $(Y^{'},\Gamma^{'}+tA_{Y^{'}}+\mathbf{M})$ has a good minimal model over $X$. Set $U_{Y^{'}}:=f^{'-1}(U)$. Once we show that $(U_{Y^{'}}+(\Gamma^{'}+tA_{Y^{'}})|_{U_{Y^{'}}}+\mathbf{M}|_{U_{Y^{'}}})$ has a good minimal model over $U$, it would follow from \cite[Theorem 1.3]{LX} that $(Y^{'},\Gamma^{'}+tA_{Y^{'}}+\mathbf{M})$ has a good minimal model over $X$. Since $\phi:X \dashrightarrow X^{'}$ is an isomorphism over $U$, $\phi_*^{-1}A^{'}|_U$ is $\mathbb{R}$-Cartier and letting $\psi:Y \dashrightarrow Y^{'}$ denote the MMP, we have $A_{Y^{'}}=\psi_*g^*A^{'}=f^{'*} \phi_*^{-1}A^{'}$ over $U$. Since $(Y^{'}, \Gamma^{'}+\mathbf{M})$ is a g-dlt modification of $(X,B+\mathbf{M})$, we have $K_{Y^{'}}+\Gamma^{'}+\mathbf{M}_{Y^{'}}=f^{'*}(K_X+B+\mathbf{M}_X)$. Combining these, we get 
        
        \begin{center}
            
        $K_{U_{Y^{'}}}+(\Gamma^{'}+tA_{Y^{'}}+\mathbf{M}_{Y^{'}})|_{U_{Y^{'}}} \newline =f^{'}|_{U_{Y^{'}}}^{*}((K_X+B+\mathbf{M}_X)|_U +t \phi_*^{-1}A|_U)$ ($**$).
        
        \end{center}
        So the generalized pair $(U_{Y^{'}}, (\Gamma^{'}+tA_{Y^{'}}+\mathbf{M}_{Y^{'}})|_{U_{Y^{'}}})$ is its own good minimal model over $U$. We conclude that $(Y^{'},\Gamma^{'}+tA_{Y^{'}}+\mathbf{M}_{Y^{'}})$ has a good minimal model over $X$. Thus it also has a log canonical model over $X$.\\

        Let $\theta:(Y^{'}, \Gamma^{'}+tA_{Y^{'}}+\mathbf{M}_{Y^{'}}) \dashrightarrow (\Tilde{X}, \Delta_{\Tilde{X}}+tA_{\Tilde{X}}+\mathbf{M}_{\Tilde{X}})$ be the birational map over $X$ to the log canonical model, where $\Delta_{\Tilde{X}}$ and $A_{\Tilde{X}}$ are the birational transforms of $\Gamma^{'}$ and $A_{Y^{'}}$ respectively. Let $\Tilde{\psi}: \Tilde{X} \rightarrow X$ be the induced birational morphism. Because of ($**$), we have $(K_{\Tilde{X}}+\Delta_{\Tilde{X}}+tA_{\Tilde{X}}+\mathbf{M}_{\Tilde{X}})|_{\Tilde{\psi}^{-1}(U)} \equiv 0$. Thus $\Tilde{\psi}:\Tilde{X} \rightarrow X$ is an isomorphism over $U$ (Note that this property wasn't enjoyed by $f^{'}: Y^{'} \rightarrow X$). If $E$ is an $f^{'}$-exceptional divisor on $Y^{'}$, then $E \cap f^{'-1}(U)\neq \emptyset$ (these are glc places), thus $E$ is contracted by $Y^{'}\dashrightarrow \Tilde{X}$ (since $\Tilde{X}$ is isomorphic to $X$ over $U$). Also note that there can't be any $\theta^{-1}$-exceptional divisors hiding in $\Tilde{X} \setminus \Tilde{\psi}^{-1}(U)$: $\theta$ is represented as $ Y^{'} \dashrightarrow \Tilde{Y} \rightarrow\Tilde{X}$, where $\Tilde{Y}$ is the corresponding good minimal model. $Y^{'} \dashrightarrow \Tilde{Y}$ is a contraction since it is a run of an actual MMP. $\Tilde{Y} \rightarrow \Tilde{X}$ is a birational morphism, hence automatically a contraction. Thus $\theta^{-1}$ can't have any exceptional divisors. We conclude that $\Tilde{\psi}$ is small.
        
        \end{proof}
\end{prop}
Now we can prove the main result of this article.\vspace{5mm}
\begin{theorem}\label{flop} \cite[Theorem 3.4]{Has} Suppose $(X,B+\mathbf{M})/S$ and $(X^{'}, B^{'}+\mathbf{M})/S$ are two generalized log canonical pairs such that $K_X+B+\mathbf{M}_X$ and $K_{X^{'}}+B^{'}+\mathbf{M}_{X^{'}}$ are nef over $S$, $\mathbf{M}_X$ and $\mathbf{M}_{X^{'}}$ are $\mathbb{R}$-Cartier and there exists a small birational map $\phi:X \dashrightarrow X^{'}$ over $S$ such that
    
    \begin{itemize}
        \item $B^{'}=\phi_*B$ and $\mathbf{M}_{X^{'}}= \phi_*\mathbf{M}_X$,
        \item there exists $ U\subset X$ open such that $\phi|_U$ is an isomorphism and all glc centers of $(X,B+\mathbf{M})$ intersect $U$
    \end{itemize}
    then (possibly after exchanging $X$ and $X^{'}$), there exist small birational morphisms from normal quasi-projective varieties $(\Tilde{X}, \Tilde{B}+\mathbf{M}) \xrightarrow{\Tilde{\psi}}(X,B+\mathbf{M})$ and $(\Tilde{X^{'}}, \Tilde{B^{'}}+\mathbf{M}) \xrightarrow{\Tilde{\psi^{'}}}(X^{'}, B^{'}+\mathbf{M})$ such that the induced birational map $(\Tilde{X}, \Tilde{B}+\mathbf{M}) \dashrightarrow  (\Tilde{X^{'}}, \Tilde{B^{'}}+\mathbf{M}) $ can be written as a composition of a finite sequence of symmetric flops over $S$ with respect to $K_{\Tilde{X}}+\Tilde{B}+\mathbf{M}$. 

    \begin{proof}
         Let $\hat{X}$ be a smooth resolution of indeterminacy of $\phi$ that extracts all glc places of $(X,B+\mathbf{M})$ and $(X^{'},B^{'}+\mathbf{M})$. Since $a(P,X,B+\mathbf{M})=a(P,X{'},B^{'}+\mathbf{M})$ for any prime divisor $P$ on $\Tilde{X}$ by negativity lemma (as in the proof of lemma \ref{mm}), $\phi(U)$ intersects all glc centers of $(X^{'},B^{'}+\mathbf{M})$. We perform induction on $\rho(\hat{X})- \max \{\rho(X), \rho(X^{'})\}$. By \cite[Lemma 3.2]{Has}, $\rho(\hat{X}) \geq \rho(X)$ and $\rho(\hat{X}) \geq \rho(X^{'})$. Thus the above quantity is atleast zero. We can assume $\rho(X) \geq \rho(X^{'})$ (by switching $X$ and $X^{'}$ if needed).\\
        
        By Proposition \ref{tech}, there exists a small projective morphism $\psi:\Tilde{X} \rightarrow X$ from a normal quasi-projective variety such that 
        $\psi$ is an isomorphism over $U$ and there is an ample $\mathbb{R}$-divisor $A^{'} \geq 0$ on $X^{'}$ such that $(X^{'},B^{'}+A^{'}+\mathbf{M})$ is glc and if $\Tilde{A}$ is the birational transform of $A^{'}$ on $\Tilde{X}$, then $K_{\Tilde{X}}+\psi_*^{-1}B+\Tilde{A}+\mathbf{M}$ is $\mathbb{R}$-Cartier and $(\Tilde{X}, \psi_*^{-1}B+\Tilde{A}+\mathbf{M})$ is glc. Since $\psi$ is small, $\mathbf{M}_{\Tilde{X}}= \psi^*\mathbf{M}_X$ is also $\mathbb{R}$-Cartier. Let $\Tilde{B}:=\psi_*^{-1}B$. Since $K_{\Tilde{X}}+\Tilde{B}+\mathbf{M}_{\Tilde{X}}=\psi^*(K_X+B+\mathbf{M})$, $(\Tilde{X},\Tilde{B}+\mathbf{M})$ is glc and $K_{\Tilde{X}}+\Tilde{B}+\mathbf{M}_{\Tilde{X}}$ is nef. Also note that $\psi^{-1}(U)$ intersects all the glc centers of $(\Tilde{X},\Tilde{B}+\mathbf{M})$ and the composition $\Tilde{X} \dashrightarrow X^{'}$ is an isomorphism on $\psi^{-1}(U)$. Thus the glc pairs $(\Tilde{X}, \Tilde{B}+\mathbf{M})$ and $(X^{'},B^{'}+\mathbf{M})$ and $\Tilde{X} \dashrightarrow X^{'}$ satisfy the hypothesis of the theorem.It is clear that $\rho(\Tilde{X}/S)- \max \{\rho(\Tilde {X} /S), \rho (X^{'}/S)\} \leq \rho(\Tilde{X}/S)-\max \{\rho(X/S), \rho(X^{'}/S)\}$. We may thus replace $(X,B+\mathbf{M})$ with $(\Tilde{X},\Tilde{B}+\mathbf{M})$. Thus we have an ample $\mathbb{R}$-divisor $A^{'} \geq 0$ on $X^{'}$ such that $(X^{'},B{'}+A^{'}+\mathbf{M})$ is glc, $A:= \phi_*^{-1}A^{'}$ is $\mathbb{R}$-Cartier and $(X,B+A+\mathbf{M})$ is glc. Note that $(X^{'},B^{'}+tA^{'}+\mathbf{M})$ is a glc model of $(X,B+tA+ \mathbf{M})$ for all $t \in (0,1]$. \\

        Suppose $K_X+B+tA+\mathbf{M}_X$ is nef for some $t \in (0,1]$. Then since $(X,B+tA+\mathbf{M})$ has a glc model, it also has a good minimal model. So by lemma \ref{mm}, any minimal model has to be good. Thus $K_X+B+tA+\mathbf{M}_X$ is semiample. Since $X$ and $X^{'}$ are isomorphic in codimension $1$, Proj $ R(K_X+B+tA+\mathbf{M}_X) \cong$  Proj $R(K_{X^{'}}+B^{'}+tA^{'}+\mathbf{M}_{X^{'}}) \cong X^{'}$. Thus $\phi:X \rightarrow X^{'}$ is a morphism given by $|K_X+B+tA+\mathbf{M}_X|_{\mathbb{R}}$ and this corresponds to the small morphism $\Tilde{\psi^{'}}$ listed in the theorem. \\

        Otherwise, if $K_X+B+tA+\mathbf{M}_X$ is not nef for all $t \in (0,1]$, we make the following claim (see \cite[Lemma 2]{Ka}): \\
        
            \emph{Claim}: In the above situation, there exists $t \in (0,1]$ such that the birational transform of $K_X+B+\mathbf{M}$ is trivial over all extremal contractions of any sequence of steps of the $(K_X+B+tA+\mathbf{M})$-MMP over $S$ \cite[Theorem 5.3]{HL}
            \begin{center}
                $(X,B+tA+ \mathbf{M})=(X_0,B_0+tA_0+\mathbf{M}) \dashrightarrow (X_1, B_1+tA_1+\mathbf{M}) \dashrightarrow \cdots \dashrightarrow (X_l,B_l+tA_l+ \mathbf{M}) $
            \end{center}

               \emph{Proof of claim:} Indeed, choose $k \in \mathbb{N}$ such that $k(K_X+B+\mathbf{M}_X)$ is Cartier. Let $e:=\dfrac{1}{2k \dim X+1}$. Since $K_X+B+etA+\mathbf{M}_X$ can't be nef, take any extremal ray $R$ negative with respect to it. Since $K_X+B+\mathbf{M}_X$ is nef, this forces $(A \cdot R) <0$. This implies $(K_X+B+tA+\mathbf{M}_X \cdot R)<0$. Since $K_X+B+tA+\mathbf{M}_X$ is glc, by \cite[Theorem 5.1]{HL}, $R$ is generated by a rational curve $C$ such that $0 >((K_X+B+tA+\mathbf{M}_X) \cdot C) \geq -2 \dim X$. We claim that $((K_X+B+\mathbf{M}_X) \cdot C)=0$. Otherwise since $k(K_X+B+\mathbf{M}_X)$ is nef and Cartier, we have $((K_X+B+\mathbf{M}_X) \cdot C) \geq \dfrac{1}{k}$. With this,

                \begin{center}
                    $(K_X+B+etA+\mathbf{M}_X) \cdot C) \newline = \dfrac{1}{2k \dim X+1} ((K_X+B+tA+ \mathbf{M}_X) \cdot C) \newline + \dfrac{2k \dim X}{2k \dim X +1}((K_X+B+\mathbf{M}_X) \cdot C )\newline \geq \dfrac{1}{2k \dim X+1} (-2 \dim X +2 \dim X) =0$,
                \end{center}
                a contradiction. Note that $k(K_{X_i}+ B_i+ \mathbf{M}_{X_i})$ is nef over $S$ and Cartier (by the basepoint free theorem) at every step of the above MMP. So the above arguments are valid if we replace $X$ with any $X_i$, $0 \leq i \leq l$. Thus we are done.\\
        
        Thus we can find $t \in (0,1]$ such that there exists a $(K_X+B+t A +\mathbf{M})$-MMP of flips over $S$

        \begin{tikzcd}
            (X,B+tA+\mathbf{M})\arrow[r,"\phi_1",dashed]&\cdots \vspace{2mm}
            \arrow[r,"\phi_l",dashed]& (X_l,B_l+tA_l+\mathbf{M})
        \end{tikzcd}

        which is trivial with respect to $K_X+B+\mathbf{M}$ and terminating with a good minimal model for $K_X+B+tA+\mathbf{M}$. Let $\psi^{'}: X_l \rightarrow X^{'}$ be the induced small birational morphism from the good minimal model to the glc model and $B^{'}:=\psi_*B_l$. Since $\phi_i$ is small for all $i$, $N^1(X_i/S)_{\mathbb{R}}$ injects into $N^1(X_{i+1}/S)_{\mathbb{R}}$ as in the proof of lemma \ref{sym}. Thus
        \begin{center}
            $\rho(X/S) \leq \rho(X_1/S) \leq \cdots \leq \rho(X_l/S)$
        \end{center}
        It remains to show that $X$ and $X_l$ are connected by symmetric flops.\\

        \underline{Case 1}: $\rho(X/S) < \rho(X_l/S)$.\vspace{3mm}
        
        Since the relative stable base locus $B(K_X+B+tA+\mathbf{M}_X/S) \subset X \setminus U$, it follows that the induced birational map $X \dashrightarrow X_l $ is an isomorphism on $U$. The image of $U$ in $X_l$ is $\psi^{'-1}(\phi(U))$. Since $K_{X_l}+B_l+\mathbf{M}_{X_l}= \psi^{'*}(K_{X^{'}}+B^{'}+\mathbf{M}_{X^{'}})$, all glc centers of $(X_l,B_l+\mathbf{M})$ intersect $\psi^{'-1}(\phi(U))$. Thus the induced birational map $(X,B+\mathbf{M}) \dashrightarrow (X_l,B_l+ \mathbf{M})$ satisfies the hypotheses of the theorem. Since $\rho(X/S) \geq \rho(X^{'}/S)$ and $\rho(X_l/S) >\rho(X/S)$, it follows that 
        \begin{center}
            $\rho(\hat{X}/S)-\max\{\rho(X/S),\rho(X^{'}/S)\} > \rho(\hat{X}/S) -\max \{\rho(X/S), \rho(X_l/S)\}$,
        \end{center} where $\hat{X}$ is the common resolution of $X$ and $X^{'}$ considered at the beginning of the proof.

        We can replace $(X^{'}, B^{'}+\mathbf{M})$ by $(X_l,B_l+\mathbf{M})$. By induction hypothesis, there exist small birational models of $X$ and $X_l$ which are connected by a sequence of symmetric flops. Thus we are done in this case.\\

        \underline{Case 2}: $\rho(X/S) = \rho(X_l/S)$.\\

        Note that $\mathbf{M}_{X_i}$ are $\mathbb{R}$-Cartier for all $i$ (since this property is preserved by flips). Thus if we look at the $i$-th flip over $S$
        
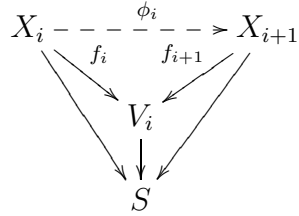
        \begin{equation*}
\xymatrix@R=16pt
{
X_i\ar@{-->}^{\phi_i}[rr]\ar[dr]^{f_i}\ar[ddr]&&X_{i+1} \ar[dl]_{f_{i+1}}\ar[ddl]\\
&V_i\ar[d]\\
&S
}
\end{equation*}

        then since $\rho(X_i/S)= \rho(X_{i+1}/S)$ and $-(K_{X_i}+B_i+tA_i+\mathbf{M}_{X_i})$ is ample over $V_i$, by lemma \ref{sym} , $\phi_i$ is a symmetric flop with respect to $K_{X_i}+B_i+\mathbf{M}_{X_i}$ for all $i$.

        \end{proof}

\end{theorem}

\section{Applications}
\begin{lemma}\label{small} \cite[Lemma 3.3]{Has}
    Let $(X, \Delta+\M)$ be a glc pair with $\M_X$ $\mathbb{R}$-Cartier and $\psi: X^{'} \to X$ a small birational morphism. Suppose there exists an open subset $U \subset X$ such that all glc centers of $(X, \Delta+\M)$ intersect $U$. Then \begin{enumerate}
        \item $R^p \psi_*\mathcal{O}_{X^{'}}=0$ for all $p>0$,
        \item for any Cartier divisor $D^{'}$ on $X^{'}$, if $D^{'} \equiv_X 0$, then $D^{'}\sim \psi^*D$ for some Cartier divisor $D$ on $X$.
    \end{enumerate}
    \begin{proof}
        We may assume codim $(X^{'} \setminus \psi^{-1}U) \geq 2$. We first show that there exists an $\mathbb{R}$-Cartier divisor $G^{'} \geq 0$ on $X^{'}$ such that $-G^{'}$ is ample over $X$ and $(X^{'}, \psi_*^{-1}\Delta+G^{'}+\M)$ is glc. Since $K_{X^{'}}+ \psi_*^{-1}\Delta +\M = \psi^*(K_X+\Delta+\M)$, $(X^{'}, \psi_*^{-1}\Delta+\M)$ is glc and all its glc centers intersect $\psi^{-1}(U)$. Pick an ample divisor $A^{'}$ on $X^{'}$ and an  ample divisor $H$ on $X$. Since $\psi$ is an isomorphism over $U$, there exists $s >0$ such that $(s \psi^*H-A^{'})|_{\psi^{-1}(U)}$ is ample. Since all glc centers of $(X^{'}, \psi_*^{-1}\Delta+\M)$ intersect $\psi^{-1}(U)$, by taking the closure of a general member of $|(s\psi^*H-A^{'})|_{\psi^{-1}(U)}|_{\mathbb{R}}$, we get $0 \leq H^{'} \sim s \psi^*H -A^{'}$ such that Supp $H^{'}$ contains no glc centers of $(X^{'}, \psi_*^{-1}\Delta+\M)$. This shows that there exists some $t>0$ such that $(X^{'}, \psi_*^{-1}\Delta+tH^{'}+\M)$ is glc and all its glc centers intersect $\psi^{-1}(U)$. Take $G^{'}:= tH^{'}$. Since $\psi $ is small, $\M_{X^{'}}= \psi^*\M_X$ is $\mathbb{R}$-Cartier and then $(X^{'}, \psi_*^{-1}\Delta+G^{'})$ is lc by Lemma \ref{M}. Since $-(K_{X^{'}}+ \psi_*^{-1}\Delta +G^{'})$ is ample over $X$, by Kodaira vanishing for lc pairs \cite[Theorem 5.6.4]{Fuj}, it follows that $R^p\psi_*\mathcal{O}_{X^{'}}=0$ for all $p>0$. \\

        Now we prove the second assertion. If $\rho(X^{'}/X) =0$, then $\psi$ is an isomorphism and there is nothing to prove. Suppose $\rho(X^{'}/X) >0$. As observed above, $(X^{'}, \psi_*^{-1}\Delta+G^{'}+\M)$ is lc and $K_{X^{'}}+\psi_*^{-1}\Delta+G^{'}+\M$ is in particular not nef over $X$. Then by contraction theorem for glc pairs \cite[Theorem 1.3 (4)]{HL}, there exists an extremal contraction $f: X^{'} \to X^{''}$ over $X$ with the property that for any Cartier divisor $D^{'}$ on $X^{'}$, if $D^{'} \equiv _X 0$, then $D^{'} \sim f^*D^{''}$ for some Cartier divisor $D^{''}$ on $X^{''}$ and $D^{''} \equiv _X 0$. The induced morphism $g: X^{''} \to X $ is small and all the glc centers of $(X^{''}, g_*^{-1}\Delta +G^{''}+\M)$ (where $G^{''}= f_*G^{'}$) intersect $g^{-1}(U)$. Since $\rho(X^{''}/X) < \rho(X^{'}/X)$, by induction hypothesis to $g$, $D^{''} \sim g^*D$ for some Cartier divisor $D$ on $X$. Then $D^{'} \sim f^*D^{''} \sim f^*g^*D = \psi^*D$.
    \end{proof}
\end{lemma}
We now have the following consequence of Theorem \ref{flop}:

\begin{corollary}\label{app} \cite[Theorem 1.2]{Has}
    Suppose $(X,B+\mathbf{M})/S$ and $(X^{'}, B^{'}+\mathbf{M})/S$ are two generalized log canonical pairs with structure morphisms $\pi: X \to S$ and $\pi^{'}: X^{'}\to S$ and such that $K_X+B+\mathbf{M}_X$ and $K_{X^{'}}+B^{'}+\mathbf{M}_{X^{'}}$ are nef over $S$, $\mathbf{M}_X$ and $\mathbf{M}_{X^{'}}$ are $\mathbb{R}$-Cartier and there exists a small birational map $\phi:X \dashrightarrow X^{'}$ over $S$ such that
    
    \begin{itemize}
        \item $B^{'}=\phi_*B$ and $\mathbf{M}_{X^{'}}= \phi_*\mathbf{M}_X$,
        \item there exists $ U\subset X$ open such that $\phi|_U$ is an isomorphism and all glc centers of $(X,B+\mathbf{M})$ intersect $U$.
    \end{itemize}
    Then we have the following:

    \begin{enumerate}
        \item $R^p \pi_* \mathcal{O}_X \cong R^p \pi^{'}_*\mathcal{O}_{X^{'}}$ for all $p >0$. In particular, if $S$ is a point, then $H^i(X, \mathcal{O}_X) \cong H^i(X^{'}, \mathcal{O}_{X^{'}})$ for all $i>0$,
        \item $K_X+B+\M_X$ and $K_{X^{'}}+B^{'}+\M_{X^{'}}$ have the same Cartier index.
        \end{enumerate}

    \begin{proof}
        By Theorem \ref{flop}, there exist small birational morphisms $f:\Tilde{X} \to X$ and $\Tilde{X^{'}}\to X^{'}$ such that $\Tilde{X}$ and $\Tilde{X^{'}}$ are connected by a sequence of symmetric flops. Let

        \begin{equation*}
\xymatrix@R=16pt
{
\Tilde{X}_i\ar@{-->}^{\phi_i}[rr]\ar[dr]^{f_i}\ar[ddr]&&\Tilde{X}_{i+1}
\ar[dl]_{f_{i+1}}\ar[ddl]\\
&V_i\ar[d]^{\pi_i}\\
&S
}
\end{equation*}
denote the $i$-th link in the flop chain. Note that $\M_{\Tilde{X_i}}$ is $\mathbb{R}$-Cartier for all $i$. By Lemma \ref{small}, $R^pf_*\mathcal{O}_X =0= R^pf^{'}_*\mathcal{O}_{\Tilde{X}^{'}}$ for all $p>0$ which, by Grothendieck spectral sequence, gives isomorphisms $R^p \pi_*\mathcal{O}_X \cong R^p(\pi \circ f)_*\mathcal{O}_{\Tilde{X}}$ and $R^p\pi^{'}_*\mathcal{O}_{X^{'}} \cong R^p (\pi^{'}\circ f^{'})_*\mathcal{O}_{\Tilde{X}^{'}}$ for all $p \geq 0$. Hence it is enough to show that $R^p(\pi_i \circ f_i)_*\mathcal{O}_{\Tilde{X}_i} \cong R^p(\pi_i \circ f_{i+1})_*\mathcal{O}_{\Tilde{X}_{i+1}}$ for all $ p \geq 0$ and $i \geq 0$. \\

Note that $f_i$ is $(K_{\X_i}+\B_i+\M)$-trivial and $(K_{\X_i}+\B_i+t\A_i +\M)$-negative (notation as in the proof of Theorem \ref{flop}). Letting $f_i=$ cont$_R$, if $\M_{\X_i} \cdot R \geq 0$, then $(K_{\X_i}+\B_i+t\A_i) \cdot R <0$. Note that by Lemma \ref{M}, $(\X_i, \B_i+t\A_i)$ is lc and so is $(\X_i, \B_i)$. Thus we can apply Kodaira vanishing for lc pairs \cite[Theorem 5.6.4]{Fuj} to get $R^pf_{i*}\mathcal{O}_{\X_i} =0$ for all $p>0$ and hence $R^p(\pi_i \circ f_i)_*\mathcal{O}_{\X_i} \cong R^p \pi_{i*}\mathcal{O}_{V_i}$ for all $p \geq 0$. If $\M_{\X_i} \cdot R <0$, then $\M_{\X_i}-\alpha \A_i \sim_{V_i}0$ for some $\alpha >0$ (by \cite[Theorem 1.2]{HL}). This gives $K_{\X_i}+\B_i \sim _{V_i} \alpha \A_i$ and thus $-(K_{\X_i}+\B_i)$ is $f_i$-ample and we can again apply Kodaira vanishing to get $R^p(\pi_i \circ f_i)_*\mathcal{O}_{\X_i} \cong R^p \pi_{i*}\mathcal{O}_{V_i}$ for all $p \geq 0$.\\

Now we argue for $f_{i+1}$. We can argue as in the proof of Lemma \ref{small} to get an effective $f_{i+1}$-anti-ample divisor on $\X_{i+1}$. Indeed, by assumption, there exists $U \subset \X_{i+1}$ open such that $\phi_{i}$ is an isomorphism over $U$ and all glc centers of $(\X_{i+1}, \B_{i+1}+\M)$ intersect $\phi_i(U)$. Then $f_{i+1}|_{\phi_i(U)}$ is an isomorphism. There exists $W \subset V_i$ large open such that $f_{i+1}^{-1}(W)$ is also large open and $f_{i+1}$ is an isomorphism over $W$. Replace $U$ with $U^{'}:=f_{i+1}^{-1}(W) \cup U$. Let $H$ be an ample divisor on $V_i$ and $A$ an ample divisor on $X_{i+1}$. Then there exists $s>0$ such that $(s f_{i+1}^*H-A)|_{U^{'}}$ is ample. Let $H^{'}_{U^{'}}$ be a general member of the linear system of $(s f_{i+1}^*H-A)|_{U^{'}}$. Letting $H^{'}$ denote its Zariski closure, $H^{'}$ does not contain any glc centers of $(\X_{i+1}, \B_{i+1}+\M)$. Thus there exists $t>0$ such that $(\X_{i+1}, \B_{i+1}+tH+\M)$ is glc. $H^{'}$ is clearly $f_{i+1}$-anti-ample. Thus we can apply Kodaira vanishing  and argue as above to get $R^p(\pi_i \circ f_{i+1})_*\mathcal{O}_{\X_{i+1}} \cong R^p \pi_{i*}\mathcal{O}_{V_i}$ for all $p \geq 0$.\\

Combining the conclusions of the above two paragraphs, we get that $R^p(\pi_i \circ f_i)_*\mathcal{O}_{\Tilde{X}_i} \cong R^p(\pi_i \circ f_{i+1})_*\mathcal{O}_{\Tilde{X}_{i+1}}$ for all $ p \geq 0$ and $i \geq 0$. As observed above, this gives our first assertion.\\

Now we prove the second assertion. Pick any Cartier divisor $D$ on $X$ such that $D \equiv _S r(K_X+B+\M)$ for some $r \in \mathbb{R}$. For $1 \leq i \leq n$, we denote the birational transforms of $D$ and $B$ on $\X_i$ by $\D_i$ and $\B_i$ respectively. Supposing $\D_i$ is Cartier and $\D_i \equiv _S r(K_{\X_i}+\B_i+\M)$, we can show that $\D_{i+1}$ is Cartier and $\D_{i+1} \equiv _S r(K_{\X_{i+1}}+\B_{i+1}+ \M_{X_{i+1}})$ as follows: $\D_i \equiv _{V_i} 0$ (since $\X$ and $\X^{'}$ are connected by a sequence of $K_{\X}+B+ \M$-flops by Theorem \ref{flop}). By contraction theorem for glc pairs \cite[Theorem 1.3(4)]{HL}, there exists a Cartier divisor $G_i$ on $V_i$ such that $\D_i \sim g_i^*G_i$. Then we have $\D_{i+1} \sim g_i^{'*}G_i$, thus $\D_{i+1}$ is Cartier and $\D_{i+1} \equiv _S r(K_{\X_{i+1}}+\B_{i+1}+\M)$ as required. Now we use this observation: since $D_0 = f^*D$ is Cartier and $D_0 \equiv _S r(K_{\X_0}+\B_0+\M)$, by induction on $i$ as above, the birational transform $\D^{'}$ of $D$ on $\X^{'}$ is Cartier and $D^{'} \equiv _S r(K_{\X^{'}}+\B^{'}+\M)$. Then $\phi_*D = f^{'}_*\D^{'}$ is Cartier by Lemma \ref{small} and $\phi_*D \equiv _S r(K_{X^{'}}+B^{'}+\M)$.

    \end{proof}
        
\end{corollary}

\section{Acknowledgements} I thank Omprokash Das for suggesting this problem, for many helpful conversations and for answering my questions. I am grateful to the anomymous referee for carefully reading the paper and making very helpful comments. Partial financial support was provided by the Department of Atomic Energy, India under project no. 12-R\&D-TFR-5.01-0500 during this work.


\begin{thebibliography}{99}

\bibitem{Bi}
C.~Birkar.
\emph{On existence of log minimal models II.}
\emph{J. Reine Angew. Math.} 658 (2011), 99–113.

\bibitem{Bi2}
C.~Birkar.
\emph{Lectures on Birational Geometry.}
https://www.dpmms.cam.ac.uk/~cb496/birgeom-paris-public.pdf
\bibitem{BZ}
C.~Birkar. D.Q.~Zhang.
\emph{Effectivity of Iitaka fibrations and pluricanonical systems of polarized pairs.}
\emph{Publ. Math. Inst. Hautes Études Sci.} 123 (2016), 283–331. 

\bibitem{Co}
A.~Corti (editor).
\emph{Flips for 3-folds and 4-folds.}
\emph{Oxford Lecture Series in Mathematics and its Applications.} 35. Oxford University Press, Oxford, 2007. x+189 pp

\bibitem{Fuj}
O.~Fujino.
\emph{Foundations of the minimal model program.}
\emph{MSJ Memoirs, 35. Mathematical Society of Japan}, Tokyo, 2017. xv+289 pp.
\bibitem{HL}
C.~Hacon, J.~Liu.
\emph{Existence of flips for generalized lc pairs.}
https://arxiv.org/abs/2105.13590

\bibitem{Has}
K.~Hashizume.
\emph{Relations between two log minimal models of log canonical pairs.}
 \emph{Internat. J. Math.}
 31 (2020), no. 13, 2050103, 23 pp. 

 \bibitem{Ka}
 Y.~Kawamata.
 \emph{Flops connect minimal models.}
\emph{Publ. Res. Inst. Math. Sci.} 44 (2008), no. 2, 419–423. 
 

\bibitem{KM}
J.~Kollár, S.~Mori.
\emph{Birational geometry of algebraic varieties.}
With the collaboration of C. H. Clemens and A. Corti. Translated from the 1998 Japanese original. Cambridge Tracts in Mathematics, 134. Cambridge University Press, Cambridge, 1998. viii+254 pp. ISBN: 0-521-63277-3

\bibitem{LX}
J.~Liu, L.~Xie.
\emph{Semi-ampleness of generalized pairs}
https://arxiv.org/abs/2210.01731


\end{thebibliography}
\end{document}